\documentclass{amsproc}

\newread\epsffilein    % file to \read
\newif\ifepsffileok    % continue looking for the bounding box?
\newif\ifepsfbbfound   % success?
\newif\ifepsfverbose   % report what you're making?
\newif\ifepsfdraft     % use draft mode?
\newdimen\epsfxsize    % horizontal size after scaling
\newdimen\epsfysize    % vertical size after scaling
\newdimen\epsftsize    % horizontal size before scaling
\newdimen\epsfrsize    % vertical size before scaling
\newdimen\epsftmp      % register for arithmetic manipulation
\newdimen\pspoints     % conversion factor
\def\epsfbox#1{\global\def\epsfllx{72}\global\def\epsflly{72}%
   \global\def\epsfurx{540}\global\def\epsfury{720}%
   \def\lbracket{[}\def\testit{#1}\ifx\testit\lbracket
   \let\next=\epsfgetlitbb\else\let\next=\epsfnormal\fi\next{#1}}%
\def\epsfgetlitbb#1#2 #3 #4 #5]#6{\epsfgrab #2 #3 #4 #5 .\\%
   \epsfsetgraph{#6}}%
\def\epsfnormal#1{\epsfgetbb{#1}\epsfsetgraph{#1}}%
\def\epsfgetbb#1{%
\openin\epsffilein=#1
\ifeof\epsffilein\errmessage{I couldn't open #1, will ignore it}\else
   {\epsffileoktrue \chardef\other=12
    \def\do##1{\catcode`##1=\other}\dospecials \catcode`\ =10
    \loop
       \read\epsffilein to \epsffileline
       \ifeof\epsffilein\epsffileokfalse\else
          \expandafter\epsfaux\epsffileline:. \\%
       \fi
   \ifepsffileok\repeat
   \ifepsfbbfound\else
    \ifepsfverbose\message{No bounding box comment in #1; using defaults}\fi\fi
   }\closein\epsffilein\fi}%
\def\epsfclipoff{\def\epsfclipstring{\ifepsfdraft\space clip\fi}}%
\epsfclipoff
\def\epsfsetgraph#1{%
   \epsfrsize=\epsfury\pspoints
   \advance\epsfrsize by-\epsflly\pspoints
   \epsftsize=\epsfurx\pspoints
   \advance\epsftsize by-\epsfllx\pspoints
   \epsfxsize\epsfsize\epsftsize\epsfrsize
   \ifnum\epsfxsize=0 \ifnum\epsfysize=0
      \epsfxsize=\epsftsize \epsfysize=\epsfrsize
      \epsfrsize=0pt
     \else\epsftmp=\epsftsize \divide\epsftmp\epsfrsize
       \epsfxsize=\epsfysize \multiply\epsfxsize\epsftmp
       \multiply\epsftmp\epsfrsize \advance\epsftsize-\epsftmp
       \epsftmp=\epsfysize
       \loop \advance\epsftsize\epsftsize \divide\epsftmp 2
       \ifnum\epsftmp>0
          \ifnum\epsftsize<\epsfrsize\else
             \advance\epsftsize-\epsfrsize \advance\epsfxsize\epsftmp \fi
       \repeat
       \epsfrsize=0pt
     \fi
   \else \ifnum\epsfysize=0
     \epsftmp=\epsfrsize \divide\epsftmp\epsftsize
     \epsfysize=\epsfxsize \multiply\epsfysize\epsftmp   
     \multiply\epsftmp\epsftsize \advance\epsfrsize-\epsftmp
     \epsftmp=\epsfxsize
     \loop \advance\epsfrsize\epsfrsize \divide\epsftmp 2
     \ifnum\epsftmp>0
        \ifnum\epsfrsize<\epsftsize\else
           \advance\epsfrsize-\epsftsize \advance\epsfysize\epsftmp \fi
     \repeat
     \epsfrsize=0pt
    \else
     \epsfrsize=\epsfysize
    \fi
   \fi
   \ifepsfverbose\message{#1: width=\the\epsfxsize, height=\the\epsfysize}\fi
   \epsftmp=10\epsfxsize \divide\epsftmp\pspoints
   \vbox to\epsfysize{\vfil\hbox to\epsfxsize{%
      \ifnum\epsfrsize=0\relax
        \includegraphics{\ifepsfdraft}%
      \else
        \epsfrsize=10\epsfysize \divide\epsfrsize\pspoints
        \includegraphics{\ifepsfdraft}%
      \fi
      \hfil}}%
\global\epsfxsize=0pt\global\epsfysize=0pt}%
{\catcode`\%=12 \global\let\epsfpercent=%\global\def\epsfbblit{%BoundingBox}}%
\long\def\epsfaux#1#2:#3\\{\ifx#1\epsfpercent
   \def\testit{#2}\ifx\testit\epsfbblit
      \epsfgrab #3 . . . \\%
      \epsffileokfalse
      \global\epsfbbfoundtrue
   \fi\else\ifx#1\par\else\epsffileokfalse\fi\fi}%
\def\epsfempty{}%
\def\epsfgrab #1 #2 #3 #4 #5\\{%
\global\def\epsfllx{#1}\ifx\epsfllx\epsfempty
      \epsfgrab #2 #3 #4 #5 .\\\else
   \global\def\epsflly{#2}%
   \global\def\epsfurx{#3}\global\def\epsfury{#4}\fi}%
\def\epsfsize#1#2{\epsfxsize}
\newtheorem{sect}{}[section]

\begin{document}

\title{Alexander Numbering of Knotted Surface Diagrams{\footnote{MRCN:57Q45}}}

\author{
J. Scott Carter \and Seiichi Kamada \and Masahico Saito}
\address{
University of South Alabama \\
Mobile, AL 36688 \\ carter@mathstat.usouthal.edu
\\University of South Alabama \\
Mobile, AL 36688 \\ skamada@mathstat.usouthal.edu 
 {\footnote{permanent address: Osaka City University 
Osaka 558-8585, JAPAN \\ 
kamada@sci.osaka-cu.ac.jp}}
\\ University of South Florida
\\ Tampa, FL 33620  \\ saito@math.usf.edu
}

\begin{abstract}
A formula that relates triple points,
branch points, and their
distance from infinity is presented.
\end{abstract}
\maketitle

\section{Introduction}

\begin{figure}
\begin{center}
\mbox{
\epsfysize=1in
\epsfbox{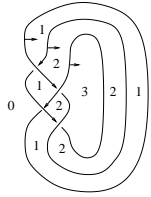}
} \end{center}
\caption{ Alexander numbering of a classical knot                
}
\label{braidLX}
\end{figure}

An {\it Alexander numbering} of a
classical knot diagram
is depicted in Fig.~\ref{braidLX}.
In this figure,
the planar regions that are separated by arcs
of the curve are labeled with integers.
The unbounded region is labeled 0.
The indices in two regions that are separated
by an arc differ by 1; the region into which
a normal arrow points has the larger index.
Such numberings were used by J.W. Alexander  \cite{Alex}
to give a combinatorial defintion of his now
well-known polynomial invariant (See also \cite{Kauf:OnKnots}).
An Alexander numbering for  
oriented 
knotted closed surfaces is 
defined similarly in Section~\ref{LXsec} 
where indices are also assigned to 
the crossing points
of a knotted surface diagram.

\begin{figure}
\begin{center}
\mbox{
\epsfxsize=3in
\epsfbox{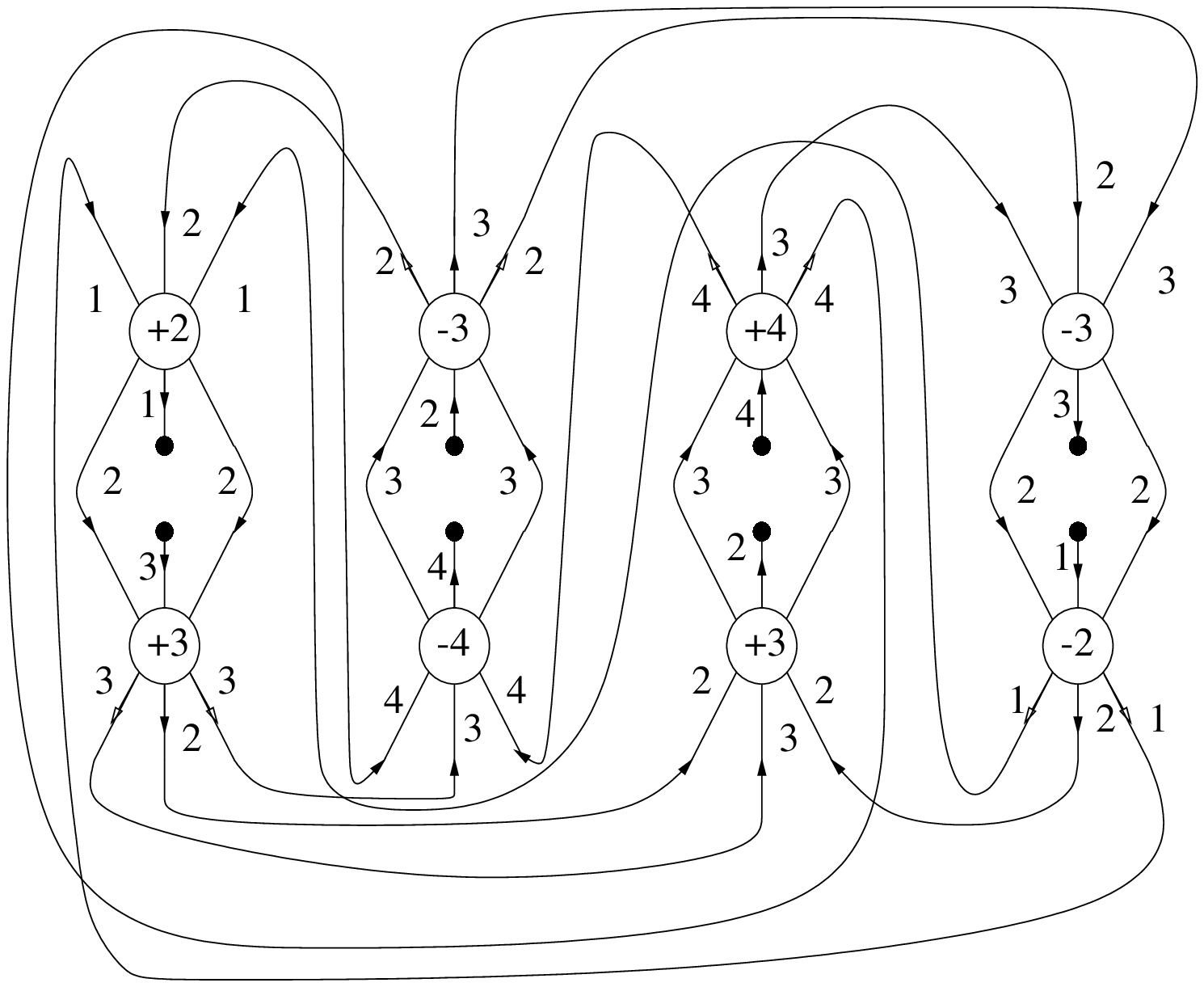}
} \end{center}
\caption{A braid chart of the 3-twist spun trefoil}               
\label{3twistchart}
\end{figure}

\begin{figure}
\begin{center}
\mbox{
\epsfxsize=2in
\epsfbox{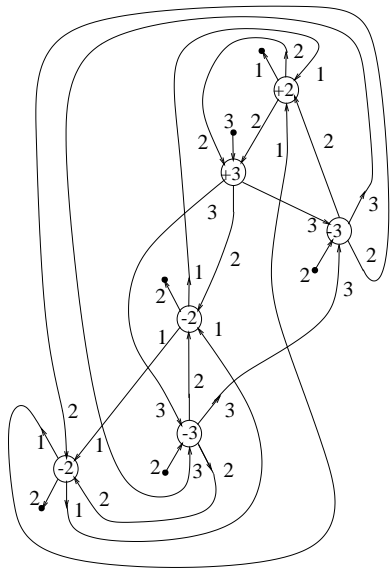}
} \end{center}
\caption{ A braid chart of Fox's \cite{FoxTrip} Example 12 } 
\label{ex12chart}
\end{figure}

 The title of this paper is an intentional pun. 
Alexander is also the given name of  JSC's second son,
 who at age 8, 
looked at the diagram of
Fig.~\ref{3twistchart} (which is SK's diagram of
the 3-twist spun trefoil and at the time was
on the blackboard of JSC's office), and observed
that the sum of the 
encircled 
labels 
added to 0. Alexander Carter's observation on this
diagram does not hold for every chart of knotted surfaces
(see Fig.~\ref{ex12chart}), but there is a relation among the
indices of the triple points and branch points that we develop herein.
These diagrams represent surface braids,
 and will be explained
in  Section~\ref{prelim}.

\begin{sect}{\sc Historical remarks.}
{\rm
A formula
 between Whitney degree \cite{Wh1}
 and
the Alexander numbering was given in \cite{MC},
which was based on the integral
calculus  with respect to the Euler characteristic \cite{Viro}.
The Alexander numbering for surfaces in $3$-space was
used in \cite{Goryu} and  \cite{MP}.
Whitney's formula
 \cite{Wh}
relates the normal Euler class to the Euler characteristic
(see also
\cite{CS1}).
He conjectured that for non-orientable surfaces the number of
achievable values
is limited.
This conjecture was proven by Massey in \cite{Massey}; see
\cite{Kam:Osaka}
for a geometric proof.
Banchoff relates the branch points to the normal
Euler class in \cite{Banchoff}.
Banchoff's triple point formula~\cite{Ban}
also relates singularities of surface maps to their intrinsic topology.
Li shows
\cite{Li}
when a daisy graph is
realizable as the multiple point set of a
general position map.
Generalizations of Banchoff's formula appear in \cite{IzMar}
and \cite{KoC}. Relation between branch points and triple points
 are found in \cite{CCS} and \cite{CS:checker}.
A further generalization is presented in this paper.
}\end{sect}

Here is how the paper will be developed.
Section~\ref{prelim} summarizes the notation also 
found in ~\cite{CS:book}. Section~\ref{LXsec} defines the method of numbering.
Section~\ref{main} contains the 
main results: Theorem~\ref{formula} and its consequences.

\begin{sect}{\sc Acknowledgements.}
{\rm
First, we thank Alexander Carter for making his helpful observation.
Seiichi Kamada is being supported by a Fellowship from the Japan 
Society for the Promotion of Science.
}\end{sect}

\begin{figure}
\begin{center}
\mbox{
\epsfysize=1in
\epsfbox{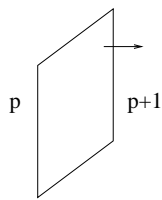}
} \end{center}
\caption{ Alexander numbering changes by 1 through a sheet               
}
\label{index}
\end{figure}

\section{Notation} 
\label{prelim}

In this section we review  knotted surface diagrams
and 
the chart description of 
surface braids 
that will be used 
in the paper. More details can be found in \cite{CS:book}.

\begin{sect} 
{\sc Generic projections and diagrams.} {\rm
For classical knots and links, under-arcs are broken in the projections
to define knot diagrams. We generalize this notion to knotted surfaces
as broken surface diagrams. 
First we develop some notation.
Let $f:F \rightarrow {\bf R}^4$ denote a smooth embedding of a closed
surface $F$ into 4-dimensional space. 

By deforming the map $f$ slightly by an ambient isotopy of ${\bf R}^4$, 
we may assume that
$f$ is a general position map. So
any point of the
image surface has a neighborhood $U$ in $3$-space such that
$p_v (f(F))   \cap U$
looks like one of the pictures indicated in 
Fig.~\ref{index} or the left of Fig.~\ref{surfaceLX}.
In other words, there is a diffeomorphism  of $U$ into ${\bf R}^3$
such that the image of $F$ coincides with the intersection of $1$,
$2$ ({\it a double point curve}), or $3$ coordinate planes
({\it a triple point}), or 
is                     
like the cone on a figure $8$ ({\it a branch point}).
The neighborhood of a double point curve, a triple point,
a branch point are depicted in Figure~\ref{surfaceLX} (A), (B), and (C),
respectively. 
When the surface is oriented, we take normal vectors $\vec{n}$
to the projection of the surface such that the triple 
$(\vec{v}_1, \vec{v}_2, \vec{n})$ matches the orientaion of 3-space,
where $(\vec{v}_1, \vec{v}_2)$ defines the orientation of the surface. 
Such normal vectors are defined on the projection at all points other than
the isolated branch points.

\begin{figure}
\begin{center}
\mbox{
\epsfxsize=2.5in
\epsfbox{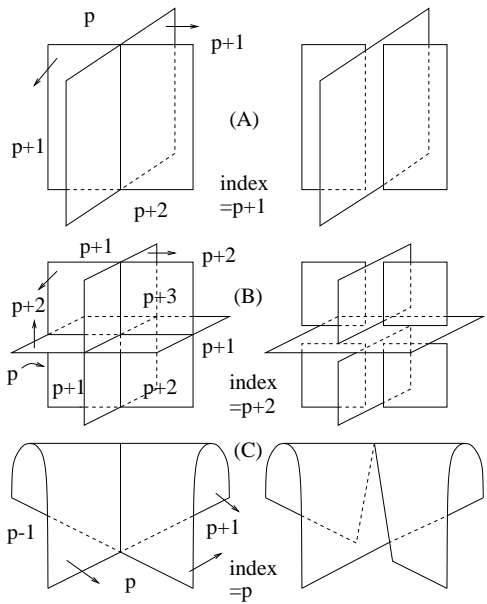}
} \end{center}
\caption{Projections, broken surface diagrams and indices }
\label{surfaceLX} 
\end{figure}

There is 
an immersion  in  ${\bf R}^3 $   
of a compact 1-manifold  with boundary 
whose image is the closure 
of the set                   
$$\{ x \in p(f(F) ) | p(f(x_1))= x = p(f(x_2))   \ \ {\mbox{\rm for some}} 
\ 
x_1 \neq x_2, \ \
\ \ {\mbox{\rm where}} \ \ x_1 , x_2 \in F \}.$$
The closure of 
the above                                    
set  is called 
{\it the double point set. } 
The image of the boundary points of the 1-manifold are the branch points 
 of the map $p\circ f.$
At a triple point three curves meet transversely, 
so a triple point is a 6-valent vertex. 

In the {\it diagram of a knotted surface}, 
one of the two sheets
along the double point arcs
to indicate that 
it is 
higher than the other sheet in the $w$-direction.
Broken surfaces are indicated in Fig.~\ref{surfaceLX}.

The double point curves are oriented in such a way that 
the orientation $\vec{v}$ together with normal vectors of top and 
bottom sheets $( \vec{n_1}, \vec{n_2}, \vec{v})$ 
matches 
the orientation of the space.
The double point curve depicted in Fig.~\ref{surfaceLX} (A) is
 oriented downward.

} \end{sect}

\begin{sect} {\sc Signs of branch and triple points.\/} {\rm
The {\it sign of a triple point} 
 is defined as follows \cite{CS:checker}.
For the normal vectors $v_1, v_2, v_3$ of top, middle, bottom 
sheet, respectively, 
if the triple $(v_1, v_2, v_3)$ matches the orientation of ${\bf R}^3$,
then 
the sign is positive,  and negative if otherwise.

	The {\it sign of  a branch point}
is defined as follows \cite{CCS,CS1}. 
If the 
the double curve ending at the given branch point is 
oriented towards (resp. 
away from) the branch point, then 
the sign is 
negative (resp. positive).
The triple point depicted in Fig.~\ref{surfaceLX} is negative, 
and the branch point is positive.

} \end{sect}

\begin{sect} {\sc  Charts of surface braids.\/}
{\rm 
The notion of surface braids, a generalization
of Artin's braid theory, was  originally proposed by Viro, and a similar 
notion had been used by Rudolph \cite{Rud}. 
We review the chart description of surface braids as developed in 
\cite{CS:bm,CS:beyond,Kam:ribbon,Kam:Nato}.

Figure~\ref{ebcproj} shows how to express certain generic surfaces
in $3$-space by means of planar graphs. 
We consider a surface $S$ 
(this $S$ corresponds to $f(F)$ in the preceding section) in a box
$B= I_1 \times I_2 \times I_3 \subset {\bf R}^3$
a schematic of which
is depicted on the right of Figure~\ref{ebcproj}, where
$I_j$ denotes a copy of the unit interval, for $j=1,2,3$.
We require that the surface $S$ in $B$
satisfies the following conditions.

\begin{itemize}

\item
$S$ is generic.

\item 
The boundary $\partial S$  of $S$ is 
a closed trivial  braid  
contained in
$\partial (I_1 \times I_2) \times I_3$.

\item
The projection $p: I_1 \times I_2 \times I_3
\rightarrow I_1 \times I_2$ restricted
to $S$ is a branched
covering such that each branch point is 
{\it simple}. 
The preimage of a branch point, then, consists of a unique
 branch point of degree 2
 and a collection of points around which the projection
 is a local homeomorphism.

\end{itemize}

Let $D$  be the double point set of $S$.
Then $p(D) \subset I_1 \times I_2$ is a planar graph.
Let us call this an {\it un-oriented chart}.
An un-oriented chart has univalent vertices corresponding to
branch points,
$4$-valent vertices that correspond to the crossings of the
projections of double arcs,
 and $6$-valent vertices corresponding to
triple points of $S$.
A generic
intersection $S\cap I_1 \times \{ t \} \times I_3$,
consists of
intersecting strings in $ I_1 \times \{ t \} \times I_3$.
Such intersections are shown in the figure by movie strips.

When the surface $S$ is the projection of a surface braid in the 4-disk
$I_1 \times I_2 \times I_3 \times I_4$, the fourth coordinate can be 
indicated in the movies by broken arcs.
Then 
the arcs in a braid chart are oriented to agree with 
the orientation 
of the corresponding double curves. A downward pointing arrow in the chart
(with respect to the height direction in the page)
that is labeled with an integer $i$ 
corresponds to the $i\/$th braid 
generator in which the $i\/$th string passes over the 
$(i+1)\/$st string.
In this way, an oriented braid chart for surface braids is defined
\cite{Kam:Nato}.

The univalent vertices of a chart correspond to branch points and are 
called {\it black vertices}.
The $6$-valent vertices correspond to triple points 
and are called {\it white vertices}. 
The signs of triple points 
on a braid chart can be determined as follows:
if the indices of arcs at a triple point is $p$ and $p+1$
and the middle incoming arc is labeled by $p+1$, 
then the triple point (white vertex) is positive,
and otherwise negative.
\begin{figure}
\begin{center}
\mbox{
\epsfxsize=3in
\epsfbox{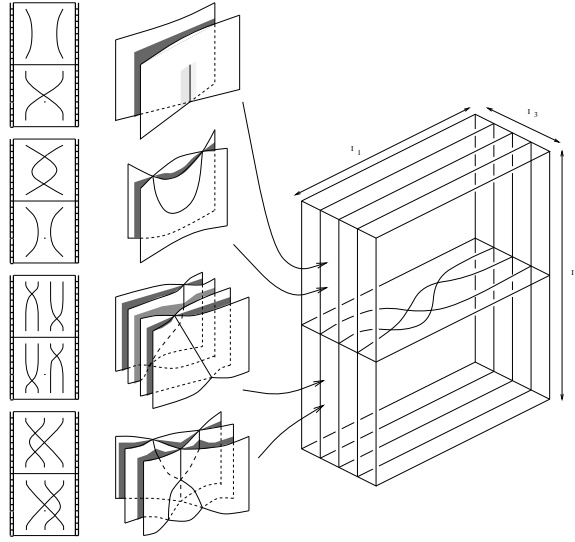}
} \end{center}
\caption{ Generic surfaces, movies,    
 and surface braids                  
}
\label{ebcproj}
\end{figure}
By capping off the nested boundary circles on the boundary of the box
$B$ by nested disjoint disks, we get a closed surface
in $3$-space. 
If it is a projection of a knotted surface, 
taking such a closure gives rise to a closed 
surface braid 
in $4$-space.

The figures~\ref{3twistchart} and \ref{ex12chart} represent 
the braid charts of 3-twist-spun trefoil (see \cite{Rolf} for example)
and Fox's Example 12 in \cite{FoxTrip}. 

} \end{sect}

\section{Alexander numberings } 
\label{LXsec}

Let $F$ denote an oriented surface that is embedded in $4$-space via a
map $f:F\rightarrow{\bf R}^4$. Consider a general position projection of
the knotting $f$ into $3$-space and the associated knotted surface
diagram.
 Choose a normal orientation in 3-space for the projection such
that the tangent orientation followed the normal orientation agrees with
the right-handed orientation of 3-space. Number with integers the
3-dimensional regions in the complement of the projection according to
the convention that (a) the unbounded region is numbered 0,
(b) regions that are separated by a 2-dimensional face are numbered
consecutively, and (c)
the normal vector to the surface points towards the region with largest
number. Such an indexing is called an {\it Alexander numbering}.
This was defined in \cite{Goryu,MP} and was called {\it index} 
or {\it degree}.

If the regions adjacent to a double arc are numbered $p$, $p+1$, $p+1$,
and $p+2$, then {\it the
Alexander number of the double arc} is $p+1$ (Fig.~\ref{surfaceLX} (A)). 
In a neighborhood of a triple point, one of the adjacent regions is
numbered $p$,
three are numbered $p+1$, three are numbered $p+2$, and one is numbered
$p+3$. The 
{\it Alexander number of  the triple point} is $p+2$ 
(Fig.~\ref{surfaceLX} (B)). 
The {\it Alexander number
of a branch point} is the number of the region that is not interior to
the figure 8. Thus if the regions interior to the figure 8
are labeled $p-1$ and $p+1$, then the Alexander number of the branch
point is $p$ (Fig.~\ref{surfaceLX} (C)). 

\begin{figure}
\begin{center}
\mbox{
\epsfxsize=2in
\epsfbox{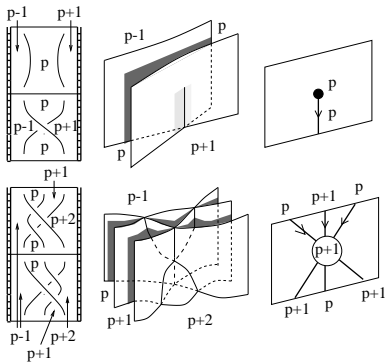}
} \end{center}
\caption{ Indices at braid generators               
}
\label{ebcLX}
\end{figure}

When a knotted surface is given as a surface braid, then the cross
sectional braids are closed with arcs that pass to the right as in
Fig.~\ref{braidLX}.
We may assume that the surface is oriented so that the normal direction
points
consistently to the right 
on the ``braid portion of the diagram''
(On the closure part, the normal points left: Fig.~\ref{braidLX}).
The {\it braid index of a branch point} is the index of the braid
generator that is created or distroyed at the black branch point in the
chart. The {\it braid index of a triple point}
is the largest braid index among the indices that appear adjacent to the
triple point. Thus at a white vertex in a braid chart 3 edges each are
incident with indices 
$p$ and $p+1$. The braid index  of the vertex
is $p+1$. A white vertex in the chart corresponds to the braid relation
$\sigma_p \sigma_{p+1} \sigma_p = \sigma_{p+1} \sigma_p \sigma_{p+1};$
the braid index is $(p+1)$.
The situation is depicted in Fig.~\ref{ebcLX}.

\begin{sect}{\sc Observation.\/} The Alexander numberings of the double,
branch, and triple points agree with the braid index of the
corresponding points when the surface is in braid form.  \end{sect}
 {\it Proof.}  The region containing $\infty$ is immediately to the left
of the first arc in a braid word. 

The Alexander numberings of double, branch, triple points are
indicated in Figs.~\ref{3twistchart} and \ref{ex12chart}
in their charts.

\section{Statement of Main Result}
\label{main}

\begin{sect}{\sc Definition.} {\rm Let $T(p, \delta)$ denote the number
 of triple points of 
Alexander index $p$ with sign $\delta= \pm.$ Let 
$B(p, \sigma)$ denote the number of branch points of 
Alexander index $p$ with sign $\sigma = \pm.$
}\end{sect}

\begin{sect}{\sc Theorem.}\label{formula}
For a diagram of a knotted oriented closed surface,
$$\sum_{p,\sigma} \sigma x_p B(p,\sigma) + 
\sum_{q,\delta} \delta y_q T(q,\delta) = 0$$
provided
$y_p= x_{p}-x_{p-1}.$\end{sect}
{\it Proof.} Consider an arc of double points of index $p$. 
Such an arc can start at a branch point 
of index $p$ or start at a triple point of index $p$ or index $(p+1).$
There are two double arcs of index $p$ that
start at each positive triple point of index $p$ 
and each negative triple of index $(p+1)$. 
One such arc starts at a negative triple point of index $p$,
positive triple of index $(p+1)$, or branch point of index $p$.
Thus the number of edges of index $p$ is
$$E(p)=B(p,+)+ 2T(p,+)+T(p,-) + T(p+1,+) + 2T(p+1,-)$$ 

The ending points of such arcs can be determined similiarly.
So,
$$E(p)=B(p,-)+ 2T(p,-)+T(p,+) + T(p+1,-) + 2T(p+1,+).$$ 

We obtain,
$$B(p,+)-B(p,-) = T(p+1,+)-T(p+1,-) -(T(p,+)-T(p,-)) 
.$$
We call this relation {\it equation} $(\star)$. 
Multiply through by $x_p$, and sum over all $p$ to get:

\begin{eqnarray*}
\sum_{p,\sigma} \sigma x_p B(p,\sigma) & =&
\sum_{p,\delta} \delta x_p T(p+1,\delta) - \sum_{p,\delta}
 \delta x_p T(p,\delta)\\
&= & \sum_{p,\delta} \delta (x_{p-1} - x_p) T(p,\delta)
\end{eqnarray*}
This completes the proof. 

\begin{sect}{\sc Corollary.} If $y_p =0$ for all $p$, then
$x_p = x_{p+1}$. Therefore,
$$\sum_{p,\sigma} \sigma B(p, \sigma) =0.$$
\end{sect}

So we recover the fact that oriented surfaces have trivial normal
 Euler classes
because the number of  signed branch points add up to 
the normal Euler number \cite{CS1}. 

\begin{sect}{\sc Corollary.} If $y_p=p$, and $x_1 = 1$, then
$x_p  = p(p+1)/2$, and 
$$\sum_{p,\sigma} \sigma \frac{p(p+1)}{2} B(p,\sigma) + 
\sum_{q,\delta} \delta q T(q,\delta) =0 .$$ 
\end{sect}

\begin{sect}{\sc Corollary.} If $y_p=1$, and $x_1 = 1$, then
$x_p = p$, and 
$$\sum_{p,\sigma} \sigma p B(p,\sigma) + 
\sum_{q,\delta} \delta  T(q,\delta) =0 .$$
\end{sect}

\begin{sect}{\sc Corollary.} If $B(p, \sigma)= 0$ for all $p$, then
$T(p,+)=T(p,-).$  
In particular, this is true if the projection of the surface is immersed.
\end{sect}

\begin{sect}{\sc Corollary.}

$T(p,+)-T(p,-)=$ $ \sum_{i < p} B(i,+)-B(i,-)$
$=   - \sum_{i \ge p} B(i,+)-B(i,-).$ 
\end{sect}
{\it Proof.} This is a consequence of equation $\star$. 
Observe that these sums are finite.

\begin{figure}
\begin{center}
\mbox{
\epsfxsize=5in
\epsfbox{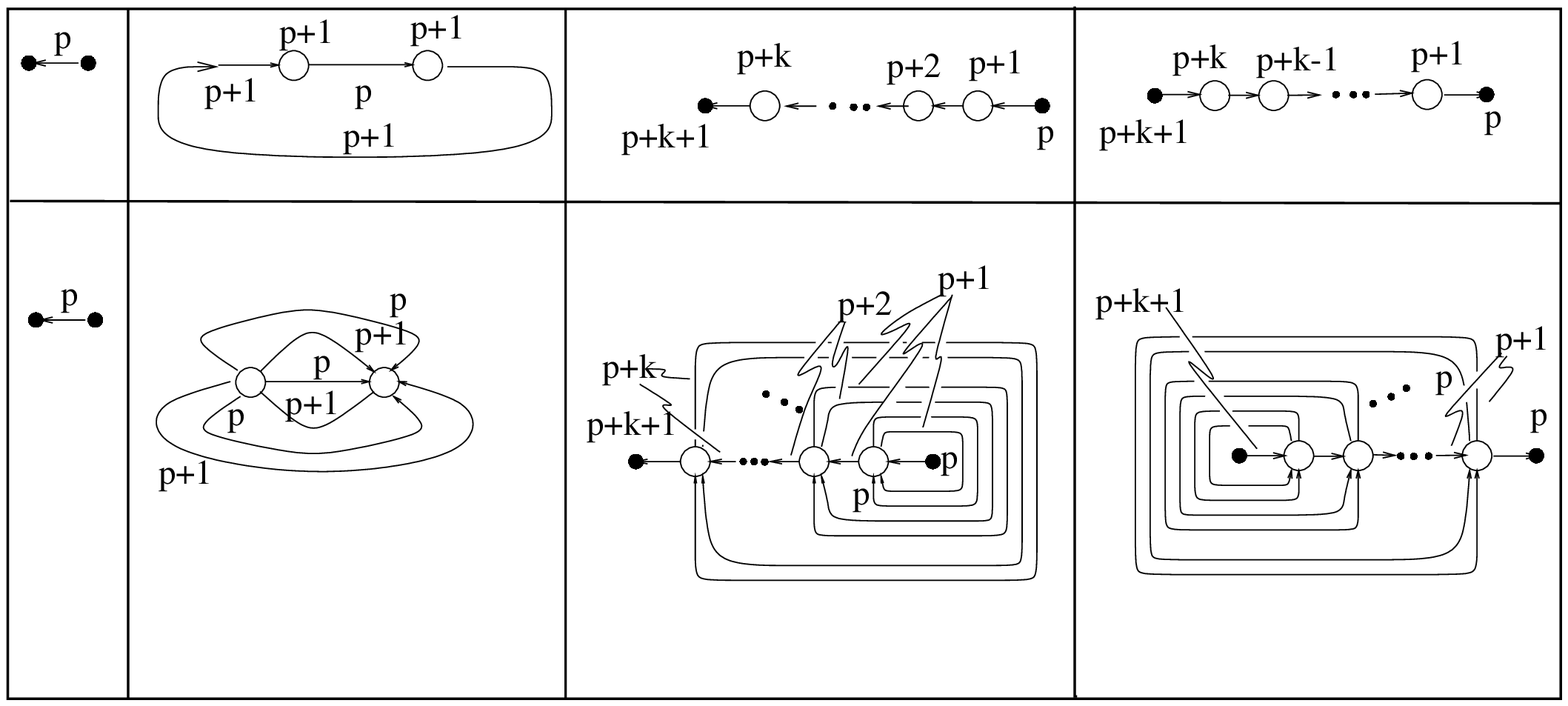}
} \end{center}
\caption{ Some charts of unknotted surfaces               
}
\label{nusomecharts}
\end{figure}

\begin{sect} {\sc Theorem.\/}
There exists a diagram of an unknot or unlink that realizes
prescibed values for $B(p,\pm)$ and $T(p,\pm)$ that satisfy the
formula in Theorem~\ref{formula}.
\end{sect} 
{\it Proof.\/}
By adding nested spheres near $\infty$, 
we may assume that the labels are all positive.
We construct a labeled planar graph with the correct
number of  black and white vertices, and then we use this 
to construct a chart of an unlink.
Between lines $x=p$ and $x=p+1$ put B(p,+) (resp. B(p,-)) black vertices 
with horizontal arrows that point away from (resp. toward) the vertex.
Put $T(p+1,+)$ (resp. $T(p+1,-)$) white vertices in the same  
region with
horizontal right-pointing (resp. left-pointing) incoming (resp. outgoing)
arrows labeled $p$ and right-pointing (resp. left pointing) outgoing (resp.
incoming) arrows labeled $(p+1)$. Along the line $x=p$ stack
$|T(p,+)-T(p,-)|$ arrows that point to the right 
if $T(p,+)-T(p,-)>0$ arrows and to the left otherwise. Connect the arrows
according to the rules that (1) lines only cross the line $x=p$ 
along the stacked arrows,  
(2) each possible pair of white vertices are connected in
a loop,  
(3) a pair of black vertices
are linearly chained by a (possibly empty) set of white vertices.  
The arcs in the graph may intersect. The only loops in the graph 
connect a pair of white vertices.

Embed the graph in the plane. 
Figure~\ref{nusomecharts} indicates how to replace the graph with bits of
charts. 
The resulting surface can be seen to be unknotted by means of chart moves.
See \cite{CS:book,Kam:ribbon} for a description of the moves.
This completes the
proof.

\begin{sect}{\sc Oriented surfaces 
with singular points.}
{\rm In this section we generalize Theorem~\ref{formula} 
to 
generically immersed oriented surfaces in $4$-space.
Such surfaces have isolated transverse double points that we call 
{\it singular points} to distinguish them 
from arcs of double points in the projection.
These can be indicated by broken surface 
diagrams as indicated in Fig.~\ref{brokdoub}.
The {\it sign of such a singular point} 
is negative if the oriented double arcs  
both point toward the singular point. The sign is positive 
if the double arcs point away from the singular point.
Say that two of the regions adjacent to the singular point have indices 
$(p+1)$, the other two have indices $p$ and $p+2$. Then the {\it index of
singular point} is $(p+1)$.

\begin{figure}
\begin{center}
\mbox{
\epsfxsize=.75in
\epsfbox{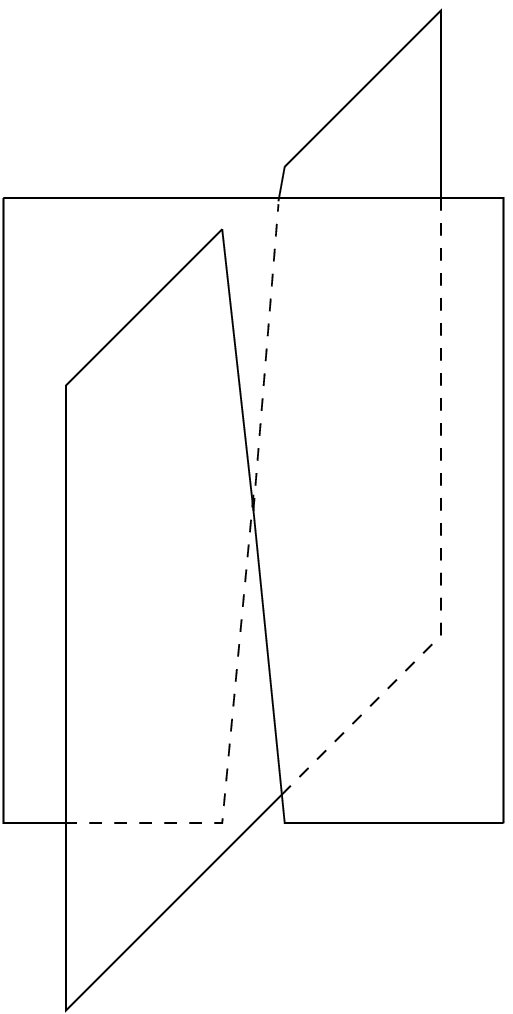}
} \end{center}
\caption{ A broken surface diagram at a singular point             
}
\label{brokdoub}
\end{figure}

Let $D(p,\epsilon)$ denote the number of singular points of 
index $p$ and sign $\epsilon$.

}\end{sect}

\begin{sect}{\sc Theorem.}
For a diagram of a knotted oriented closed surface with singular points,
\begin{eqnarray*}
\sum_{p,\sigma} \sigma x_p B(p,\sigma) + 
\sum_{q,\delta} \delta y_q T(q,\delta) +
 \sum_{r, \epsilon} \epsilon z_r D(r,\epsilon)
& =& 0 \end{eqnarray*}
provided that $y_p=x_p -x_{p-1}$ and $z_p=2 x_p.$
\end{sect}
{\it Proof.} By counting edges of label $p$, as before, we obtain
$$E(p)=B(p,+)+ 2T(p,+)+T(p,-) + T(p+1,+) + 2T(p+1,-) +2D(p,+)$$ 
$$=B(p,-)+ 2T(p,-)+T(p,+) + T(p+1,-) + 2T(p+1,+) +2 D(p,-).$$ 
Then
$$B(p,+)-B(p,-)+ 2(D(p,+) - D(p,-))$$
$$ = T(p+1,+)-T(p+1,-) -(T(p,+)-T(p,-)).$$
Multiply through by $x_p$, and sum over all $p$ to 
obtain the result.

%%%%%%%%%%%%%delete but keep for future use
%\begin{eqnarray*}
%\sum_{p,\sigma} \sigma x_p B(p,\sigma)
%+\sum_{p,\epsilon} \epsilon 2x_p D(p,\epsilon)
% & =&
%\sum_{p,\delta} \delta x_p T(p+1,\delta) - 
%\sum_{p,\delta} \delta x_p T(p,\delta)\\
%&= & \sum_{p,\delta} \delta (x_{p-1} - x_p) T(p,\delta)
%\end{eqnarray*}
%This completes the proof. 
%%%%%%%%%%%%%%%%%%%%%%%%%%%%%%%%%%%%%

\end{document}